\documentclass[11pt]{article}  
\usepackage{graphicx}
\usepackage{float,pict2e}
\usepackage[<pict2e options>]{curve2e}
\usepackage{amsmath}
\usepackage{amssymb}
\usepackage{verbatim}
\def\p{\mbox{\rm planar}}

\newtheorem{theorem}{\bf Theorem}[section]
\newtheorem{lemma}{\bf Lemma}[section]

\begin{document}

\title{On minimal triangle-free planar graphs with prescribed 1-defective chromatic number 
}


\author{Nirmala Achuthan \and N.R. Achuthan \and G. Keady.\\
Department of Mathematics and Statistics, Curtin University\\
Bentley  WA 6102,  \textsc{Australia}
}



\date{\today} 

\maketitle

\begin{abstract}
 A graph is $(m,k)$-colourable if its vertices can be coloured with
$m$ colours such that the maximum degree of the subgraph induced on
the set of all vertices receiving the same colour is at most $k$.
The $k$-defective chromatic number $\chi_k(G)$ is the least positive
integer $m$ for which graph $G$ is $(m,k)$-colourable. Let
$f(m,k;\p)$ be the smallest order of a triangle-free planar graph
such that $\chi_k(G)=m$. In this paper we show that $f(3,1;\p)=11$.
\end{abstract}


\par\noindent{\bf  Math Review Codes.}
{MSC 05C15 \and MSC 05C35}

\section{Introduction}
\label{intro} In this paper we consider  undirected graphs with no
loops or multiple edges. For all undefined concepts and terminology
we refer to~\cite{CLZ}.

If $u$ is a vertex of $G$ then $d_G(u)$, $N_G(u)$ and $N_G[u]$
denote respectively the degree, the neighbourhood, and the closed
neighbourhood of $u$ in $G$. Let $\varepsilon(G)$ denote the number
of edges in $G$.

Let $k$ be a nonnegative integer. A subset $U$ of $V(G)$ is {\it
$k$-independent} if $\Delta(G[U])\le{k}$. A $0$-independent set is
an independent set in the usual sense. A graph $G$ is {\it
$(m,k)$-colourable} if it is possible to assign $m$ colours, say
$1,2,\ldots,m$ to the vertices of $G$, one colour to each vertex,
such that the set of all vertices receiving the same colour is
$k$-independent. The smallest integer $m$ for which $G$ is
$(m,k)$-colourable is called the {\it $k$-defective chromatic
number} of $G$ and is denoted by $\chi_k(G)$. A graph $G$ is said to
be {\it $(m,k)$-critical} if $\chi_k(G)=m$ and $\chi_k(G-u) < m$ for
every $u$ in $V(G)$. A graph $G$ is said to be {\it
$(m,k)$-edge-critical} if $\chi_k(G)=m$ and $\chi_k(G-e) < m$ for
every $e$ in $E(G)$.

It is easy to see that the following statements are equivalent.
\begin{enumerate}
\item[(i)] $G$ is $(m,k)$-colourable.
\item[(ii)] There exists a partition of $V(G)$ into $m$ sets each of which is
$k$-independent.
\item[(iii)] $\chi_k(G)\le{m}$.
\end{enumerate}
Clearly $\chi_0(G)$ is the usual chromatic number. It is easy to see
that $\chi_k(G)\le\lceil\frac{|V(G)|}{k+1}\rceil$ where $|V(G)|$ is
the order of $G$. Furthermore, $\chi_k(G)\le\chi_{k-1}(G)$ for all
integers $k\ge{1}.$

The concept of $k$-defective chromatic number has been extensively
studied in the literature
(see~\cite{AAS11,CGJ,FJLS,Fr,GH,HS,SAA97b,Wo}). Given a positive
integer $m$, there exists a triangle-free graph with $G$ with
$\chi_k(G)=m$. A natural question that arises is: what is the
smallest order of a triangle-free graph $G$ with $\chi_k(G)=m$? We
denote this smallest order by $f(m,k)$. The parameter $f(m,0)$ has
been studied by several authors (see~\cite{Avis,Chvatal,GKVS,JR,HM})
and  $f(m,0)$ is determined for $m\le{5}$. It has also been shown
that $f(3,k)\le 4k+5$ for all $k\ge 0$; $f(3,1)=9$ and $f(3,2)=13$
(see ~\cite{SAA97b,AAS11}). In the same papers the corresponding
extremal graphs have been characterized. In a recent paper
~\cite{NAG12} the authors characterized triangle-free graphs on 10
vertices with $\chi_1(G)=3$. They proved that every triangle-free
graph G of order 10 with $\chi_1(G)=3$, except one, contains one of
the four (3,1)-critical triangle-free graphs of order 9.
Furthermore, the exceptional triangle-free graph on 10 vertices is a
(3,1)-edge-critical graph.

 Gr{\"{o}}tzsch~\cite{Gr} proved that if $G$ is a triangle-free
planar graph then $\chi_0(G)\le{3}$. Hence $\chi_k(G)\le{3}$ for any
triangle-free planar graph $G$ and $k\ge{1}.$

We define $f(m,k;\p)$ to be the smallest order of a triangle-free
planar graph $G$ with $\chi_k(G)=m$. Note that $f(m,k;\p)\ge f(m,k)
$. The problem of determining $f(m,k;\p)$  is relevant only for
$m=2$ and $m=3$. Considering the bipartite graph $K_{2,k}$ it is
easy to see that $f(2,k;\p)=k+2$. Clearly $f(3,0;\p)=5$.  In this
paper we prove that $f(3,1;\p)=11$.

In all the figures in this paper a double line between sets $X$ and
$Y$ means that every vertex of $X$ is adjacent to every vertex of
$Y$.
\section{Preliminary results}
\label{sec:1} We need the following results, proofs of the theorems
being in the papers cited.

\begin{theorem} (\cite{HS,Lovasz})
\label{thm:Lovasz} For any graph $G$,
$$ \chi_k(G) \le \lceil \frac{\Delta(G)+1}{k+1}\rceil
= 1 +  \lfloor\frac{\Delta(G)}{k+1}\rfloor . $$
\end{theorem}

\medskip

\begin{theorem} (\cite{Fr})
\label{thm:Frick}
If $G$ is a $(\alpha,k)-$critical graph
 then $\delta(G)\ge{\alpha-1}$.
\end{theorem}
\medskip
\begin{theorem} (\cite[7.1.11]{West})
\label{thm:West} If $l(F_i)$ denotes the length of face $F_i$ in a
plane graph $G$, then $2\varepsilon(G)= \sum_i{l(F_i)}$.
\end{theorem}

\medskip

\begin{figure}[H]
 \centering
 \setlength{\unitlength}{0.7cm}
 \begin{picture}(15,10)(0,-7.3)
 \multiput(0,0)(8,0){2}{%
  \multiput(0,0)(0,-5){2}{%
  \put(0,0){\circle{0.7}
            \makebox(0,0)[r]{$u$}}
  \multiput(1,1)(4,0){2}{%
       \polyline(0,0)(2,0)(2,1)(0,1)(0,0)
       \multiput(0.5,0.5)(1,0){2}{\circle*{0.12}}}
  \put(1,1){%
       \put(0,-3){%
            \multiput(0.5,0.5)(1,0){2}{\circle*{0.12}}}}
  \put(3,1.5){%
            \multiput(0,0.05)(0,-0.1){2}{\line(1,0){2}}}
  \put(1,1){%
       \put(0.5,0.6){\makebox(0,0)[b]{$u_1$}}
       \put(1.5,0.6){\makebox(0,0)[b]{$u_2$}}}
  \put(5,1){%
       \put(0.5,0.6){\makebox(0,0)[b]{$z_1$}}
       \put(1.5,0.6){\makebox(0,0)[b]{$z_2$}}}
  \put(1,-2){%
       \put(0.5,0.4){\makebox(0,0)[t]{$u_3$}}
       \put(1.5,0.4){\makebox(0,0)[t]{$u_4$}}}
  \polyline(.240,.254)(1,1.394)
  \polyline(.142,.320)(1,1.606)
  \polyline(.240,-.254)(1,-1.394)
  \polyline(.142,-.320)(1,-1.606)
  }}
 \multiput(0,0)(8,0){2}{%
  \put(5,-2){%
            \multiput(0.5,0.5)(1,0){2}{\circle*{0.12}}}
  \multiput(1,1)(4,0){2}{%
       \put(0,-3){%
            \polyline(0,0)(2,0)(2,1)(0,1)(0,0)}}
  \put(3,-1.5){%
            \multiput(0,0.05)(0,-0.1){2}{\line(1,0){2}}}
  }
  \put(0,0){%
    \put(5,-2){%
       \put(0.4,0.4){\makebox(0,0)[rt]{$z$}}
       \put(1.5,0.4){\makebox(0,0)[t]{$z_3$}}}
    \polyline(5.5,1.5)(5.5,-1.5)(6.5,1.5)
    \put(3.5,-2.2){\makebox(0,0)[t]{$G_1$}}}
  \put(8,0){%
    \put(5,-2){%
       \put(0.5,0.4){\makebox(0,0)[t]{$z$}}
       \put(1.5,0.4){\makebox(0,0)[t]{$z_3$}}}
    \polyline(6.5,-1.5)(5.5,1.5)(5.5,-1.5)(6.5,1.5)
    \put(3.5,-2.2){\makebox(0,0)[t]{$G_2$}}}
  \put(0,-5){%
    \put(5,-2){%
            \multiput(0.5,0.5)(1,0){2}{\circle*{0.12}}}
    \multiput(1,1)(4,0){2}{%
       \put(0,-3){%
            \polyline(0,0)(2,0)(2,1)(0,1)(0,0)}}
    \put(3,-1.5){%
            \multiput(0,0.05)(0,-0.1){2}{\line(1,0){2}}}
    \put(5,-2){%
       \put(0.5,0.4){\makebox(0,0)[t]{$z$}}
       \put(1.5,0.4){\makebox(0,0)[t]{$z_3$}}}
       \multiput(5.95,1)(0.1,0){2}{\line(0,-1){2}}
    \put(3.5,-2.2){\makebox(0,0)[t]{$G_3$}}}
  \put(8,-5){%
    \put(1,-2){%
            \polyline(0,0)(3,0)(3,1)(0,1)(0,0)}
    \put(1,-2){%
       \put(2.5,0.5){\circle*{0.12}}
       \put(2.5,0.4){\makebox(0,0)[t]{$u_5$}}}
    \multiput(0,0)(0,0.1){2}{%
       \polyline(4,-1.55)(5.654,-1.55)}
    \multiput(0,0)(0.1,0){2}{%
       \polyline(5.950,-1.154)(5.950,1)}
    \put(6,-1.5){\circle{0.7}
            \makebox(0,0)[r]{$z$}}
    \put(3.5,-2.2){\makebox(0,0)[t]{$G_4$}}}
 \end{picture}
\caption{\label{fig:one} The critical graphs of order 9 with
$\chi_1(G)=3$:
 $G_1$ to $G_4$ of~\cite{SAA97b}.
}
\end{figure}
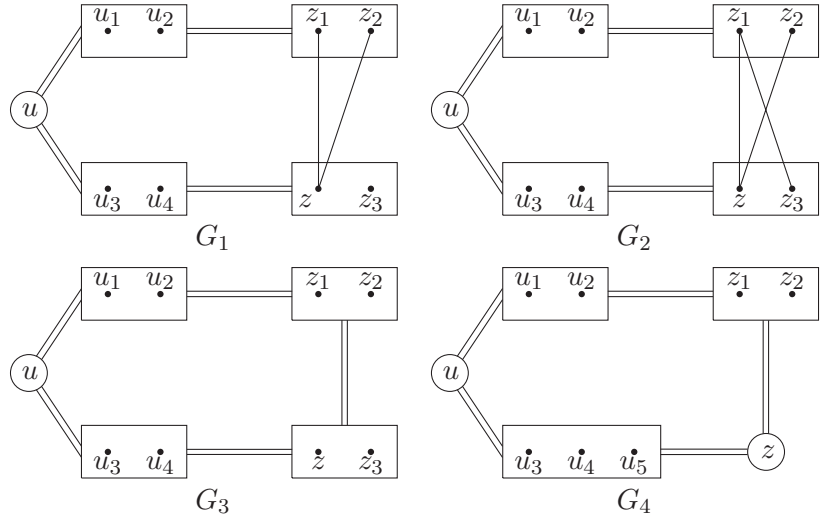

\begin{theorem} (\cite{SAA97b})
\label{thm:SAAb} The smallest order of a triangle-free graph with
$\chi_1(G)=3$ is 9, that is $f(3,1)=9$. Moreover, $G$ is a
triangle-free graph of order 9 with $\chi_1(G)=3$ if and only if it
is isomorphic to one of the graphs $G_i$, $1\le{i}\le{4}$ given in
Figure \ref{fig:one}.
\end{theorem}

\medskip
\begin{theorem}(~\cite{NAG12})
\label{thm:four} Let $G$ be a triangle-free  graph of order 10 with
$\chi_1(G)=3$. Then either $G\cong{G_5}$ (Figure~\ref{fig:g5}) or
there exists a vertex $u^*$ such that $G-u^*\cong{G_i}$ for some
$i$, $1\le{i}\le{4}$.
\end{theorem}


\begin{figure}[H]
 \centering
 \setlength{\unitlength}{0.7 cm}
 \begin{picture}(10,7)(0.0,1.5)
  \put(0.5,3){%
  \polyline(0,0)(1.5,0)(1.5,1.5)(0,1.5)(0,0)}
 \put(0.7,4.2){\makebox(0,0)[l]{$u$}}
 \put(0.7,3.3){\makebox(0,0)[l]{$v$}}
 \put(1.2,4.2){\circle*{0.2}}
 \put(1.2,3.3){\circle*{0.2}}
  \put(3.5,4){%
  \polyline(0,0)(1.5,0)(1.5,3)(0,3)(0,0)}
\put(3.7,6.5){\makebox(0,0)[l]{$u_5$}}
 \put(3.7,5.5){\makebox(0,0)[l]{$u_4$}}
 \put(3.7,4.5){\makebox(0,0)[l]{$u_3$}}
 \put(4.2,6.5){\circle*{0.2}}
 \put(4.2,5.5){\circle*{0.2}}
 \put(4.2,4.5){\circle*{0.2}}
 \put(6.5,4.5){%
  \polyline(0,0)(1.5,0)(1.5,1.5)(0,1.5)(0,0)}
\put(7.6,5.6){\makebox(0,0)[r]{$z_2$}}
 \put(7.6,4.8){\makebox(0,0)[r]{$z_1$}}
 \put(7.1,5.6){\circle*{0.2}}
 \put(7.1,4.8){\circle*{0.2}}
\put(3.5,1.0){%
  \polyline(0,0)(1.5,0)(1.5,2)(0,2)(0,0)}
\put(3.7,1.5){\makebox(0,0)[l]{$u_1$}}
 \put(3.7,2.5){\makebox(0,0)[l]{$u_2$}}
 \put(4.2,1.5){\circle*{0.2}}
 \put(4.2,2.5){\circle*{0.2}}
 \put(7.25,2.0){\circle{0.8}
            \makebox(0,0)[r]{$z$}}
 \polyline(4.2,6.5)(7.1,5.6)(4.2,4.5)(7.1,4.8)(4.2,5.5)
 \polyline(2.0,4.23)(3.5,5.53)
  \polyline(2.0,4.05)(3.5,5.35)
  \polyline(2.0,3.7)(3.5,2.2)
  \polyline(2.0,3.5)(3.5,2.0)
 \polyline(7.25,4.5)(7.25,2.4)
  \polyline(7.1,4.5)(7.1,2.36)
  \polyline(5,2.15)(6.85,2.15)
  \polyline(5,2.0)(6.84,2.0)
 \end{picture}
\caption{\label{fig:g5} $G_5$}
\end{figure}

The result of Theorem~\ref{thm:four} can also be obtained by computations,
using the methods described in the next section.

\section{Computations}

The result of Theorem~\ref{thm:six} was first obtained by computions described here.

For convenience we use the abbreviations ``tfp graph'' for a
 triangle-free planar graph and
``mtfp graph'' for a maximal tfp graph, maximal in the sense that if
another edge were to be added the graph so formed would either fail
to be triangle-free or fail to be planar or both. Recall that
$\chi_{k}(G) \le \chi_{k}(H)$ whenever $G$ is a subgraph of $H$.
Hence if $G$ is a tfp graph with $\chi_{1}(G)=3$, we can add edges
to $G$ to form a maximal tfp graph $\tilde{G}$ with
$\chi_{1}(\tilde{G})=3$.

To establish Theorem~\ref{thm:six} (and rather more)
all 11-vertex tfp graphs were generated using nauty~\cite{bdm}.
 These were then read into
a computer algebra system and tested to find which ones were not
$(2,1)-$colourable. 

\medskip

The algorithm to generate a $(2,1)-$colouring of
a triangle-free graph is provided by the following Theorem.

\goodbreak

\begin{theorem}{\label{thm:oneDefectiveColor}}
A graph $G$ of order $n$ is $(2,1)-$colourable (i.e. has
$\chi_1(G)\le{2}$) if and only if the system of  equations
\begin{eqnarray*}
b_j^2
&=& 1\qquad {\rm for\ all\ } j\in V(G)\\
b_i b_j + b_j b_k + b_k b_i
 &=&2 \qquad{\rm for\ all\ } {\rm\  paths\ } (i,j,k){\rm\ \ of\ length\ 2\ }
\end{eqnarray*}
has a solution in $Z_3^n$ where $Z_3$ is a finite field.\\
\end{theorem}

\noindent{\it Proof.} In $Z_3$, with $b_l^2=1$ for all $l=i,j,k$, we
have,
$$ b_i b_j + b_j b_k + b_k b_i =   (b_i-b_j)^2 + (b_j-b_k)^2 + (b_k-b_i)^2 .$$
\smallskip
Consider a $(2,1)-$colouring of $G$ using colours $1$ and $2$. We
define $b_i$ to be $1$ or $2$ according as the vertex $i$ is
assigned colour $1$ or $2$.

Now it is easy to see that precisely 2 of the 3 vertices on any path
$(i,j,k)$ are assigned the same colour. Thus the equations of
Theorem~\ref{thm:oneDefectiveColor} are satisfied. The converse is
easy to establish. \hfill $\square$

We remark that the $(2,1)-$colouring problem can be formulated in a
fashion similar to Theorem~\ref{thm:oneDefectiveColor} but using
systems of quadratics over $Z_2$ rather than over $Z_3$.
The connection with 3-SAT computation is given in~\cite{FJLS,GH}.

 Solving, with $n=11$,  the system of equations given in
Theorem~\ref{thm:oneDefectiveColor} with respect to each of  tfp
graphs of order 11 generated by nauty,  we found that there are
exactly six graphs with 1-defective chromatic number equal to 3.
These are presented as $G_{pi}, 1\le i \le 6$ in
Figure~\ref{fig:thmseven}. 
All these six graphs have exactly 17 edges and they
are all mtfp graphs. 

\medskip

As mentioned above, Theorem~\ref{thm:four} can be established computationally.
This time one begins by using nauty to generate all 10-vertex triangle-free graphs, not merely planar ones.
 Solving, with $n=10$,  the syste m of equations given in
Theorem~\ref{thm:oneDefectiveColor} with respect to each of  triangle-free graphs.
graphs of order 10 generated by nauty,  we establish Theorem~\ref{thm:four}.

\section{Main results}
\label{sec:2} We first prove the following useful lemmas.
\medskip
\begin{lemma}
\label{lem:one} {\it In any $(2,k)$-colouring of $K_{2,\ell}$ with
$\ell\ge{2k+1}$ the two vertices of degree $\ell$ must necessarily
receive the same colour.} Furthermore, $K_{2,\ell}$ with
$\ell\ge{3}$ is uniquely $(2,1)$-colourable.
\end{lemma}
\smallskip
\noindent{\it Proof.} Consider a $(2,k)$-colouring of $K_{2,\ell}$.
Let $C_1$ and $C_2$ be the colour classes of this ($2,k)$-colouring.
Without loss of generality let $|C_1|\ge{|C_2|}$. Clearly
$|C_1|\ge{k+2}$. If $C_1$ contains exactly one vertex of degree $l$,
then $C_1$ is not $k$-independent. Hence either $C_1$ contains both
the vertices of degree $l$ or contains neither vertex of degree $l$.
Now if $k=1$, clearly $C_2$ must contain both the vertices of degree
$l$ and $|C_2|=2$. This implies that the graph is uniquely
$(2,1)$-colourable. This establishes the lemma. \hfill $\square$
\medskip

The proof of Theorem~\ref{thm:six} uses the above result with $k=1$
and $l=3$.

\begin{lemma}
\label{lem:two}
 The graphs $G_i$, $1\le{i}\le{5}$ shown in
Figures~\ref{fig:one} and ~\ref{fig:g5} are nonplanar.
\end{lemma}
\smallskip

\begin{figure}[H]
 \centering
 \setlength{\unitlength}{1.4cm}
 \begin{picture}(2.2,3.6)(-0.2,-0.3)
 \polyline(0.,0.)(2,0.)(0.,1)(2,1)(0.,2)(2,2)(0.,0.)(2,1)
 \polyline(2,0.)(0.,2)
 \polyline(0.,1)(2,2)
 \put(-.07,-.07){\makebox(0,0)[rt]{$z_2$}}
 \put(2.07,-.07){\makebox(0,0)[lt]{$z$}}
 \put(-.10,1.00){\makebox(0,0)[r]{$z_1$}}
 \put(2.10,1.00){\makebox(0,0)[l]{$u_2$}}
 \put(-.07,2.07){\makebox(0,0)[rb]{$u$}}
 \put(2.07,2.07){\makebox(0,0)[lb]{$u_1$}}
 \put(.26,1.60){\makebox(0,0)[rt]{$u_3$}}
 \put(0.,0.){\circle*{0.1}}
 \put(2,0.){\circle*{0.1}}
 \put(0.,1){\circle*{0.1}}
 \put(2,1){\circle*{0.1}}
 \put(0.,2){\circle*{0.1}}
 \put(2,2){\circle*{0.1}}
 \put(.333,1.667){\circle*{0.1}}
 \put(1.2,-.6){\makebox(0,0)[l]{$(a)$}}
 \end{picture}
 \hspace{1.5cm}
  \begin{picture}(2.2,3.6)(-0.2,-0.3)
 \polyline(0.,0.)(2,0.)(0.,1)(2,1)(0.,2)(2,2)(0.,0.)(2,1)
 \polyline(2,0.)(0.,2)
 \polyline(0.,1)(2,2)
 \put(-.07,-.07){\makebox(0,0)[rt]{$v$}}
 \put(2.07,-.07){\makebox(0,0)[lt]{$u_3$}}
 \put(-.10,1.00){\makebox(0,0)[r]{$u$}}
 \put(2.10,1.00){\makebox(0,0)[l]{$u_2$}}
 \put(-.07,2.07){\makebox(0,0)[rb]{$z$}}
 \put(2.07,2.07){\makebox(0,0)[lb]{$u_1$}}
 \put(.26,1.60){\makebox(0,0)[rt]{$z_1$}}
 \put(0.,0.){\circle*{0.1}}
 \put(2,0.){\circle*{0.1}}
 \put(0.,1){\circle*{0.1}}
 \put(2,1){\circle*{0.1}}
 \put(0.,2){\circle*{0.1}}
 \put(2,2){\circle*{0.1}}
 \put(.333,1.667){\circle*{0.1}}
 \put(1.2,-.6){\makebox(0,0)[l]{$(b)$}}
 \end{picture}
 \caption{\label{fig:1k33one} The subdivisions of $K_{3,3}$ used in Lemma~\ref{lem:one}.}
\end{figure}
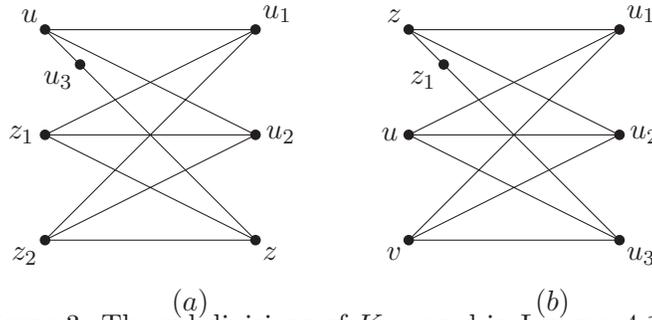

\noindent{\it Proof.} Note that the subgraph induced on the set
$\{u,u_1,u_2,u_3,z,z_1,z_2\}$ in both $G_1$ and $G_4$ is a
subdivision of $K_{3,3}$ as illustrated in $(a)$ of
Figure~\ref{fig:1k33one}. Hence $G_1$ and $G_4$ are nonplanar.


Since $G_2\cong{G_1}+(z_1,z_3)$ and $G_3\cong{G_2}+(z_2,z_3)$, $G_2$
and $G_3$ are also nonplanar.
 The subgraph of $G_5$ induced on the set
$\{u,v, u_1,u_2,u_3,z,z_1\}$ is a subdivision of $K_{3,3}$ as
illustrated in $(b)$ of Figure~\ref{fig:1k33one}. Hence the graph
$G_5$ is also nonplanar.

This proves the lemma. \hfill $\square$

\smallskip


\begin{lemma}
\label{lem:three} Let $G$ be a maximal triangle-free planar (mtfp)
graph of order $n$ with at least one odd cycle. Then the number of
edges of $G$, $\varepsilon(G)$ satisfies
\begin{equation}
\varepsilon(G)= 2n-4 -\frac{f_5}{2}\le 2 n-5 \label{eq:mtfpedge}
\end{equation}
where $f_5$ is the number of faces with 5 edges in any planar
embedding of $G$.
\end{lemma}
\goodbreak

\noindent{\it Proof.} Let $C$ be a shortest odd cycle and let $\ell$
be the length of $C$. Note the $C$ is chordless. If $\ell\ge{7}$
then we can add a chord without creating a $K_3$ or destroying the
planarity of $G$. This contradicts the maximality of $G$. Hence
$\ell=5$.

Consider a planar embedding $G'$ of $G$. Clearly $\varepsilon(G')=
\varepsilon(G)$ and every face of $G'$ is of length $4$ or $5$. Let
$f_i$ be the number of faces with $i$ edges, $i=4$ or $5$. From
\cite[7.1.13]{West} it follows that $f_5>0$.

Using Euler's formula we have
\begin{equation}
n-\varepsilon(G')+f_4 +f_5 = 2\ . \label{eq:euler}
\end{equation}
Now using Theorem ~\ref{thm:West},
\begin{equation}
4 f_4 + 5 f_5 = 2\varepsilon(G')\ . \label{eq:faceedge}
\end{equation}
Elimination of $f_4$ from equations~(\ref{eq:euler}) and
(\ref{eq:faceedge}) yields the equation at the left
of~(\ref{eq:mtfpedge}). From this it follows that $f_5$ is even.
Since $f_5>0$ we have $f_5\ge{2}$. This then yields the inequality
at the right of~(\ref{eq:mtfpedge}), and completes the proof. \hfill
$\square$

\medskip
We remark that inequality~(\ref{eq:mtfpedge}) remains valid for any
tfp graph of order $n$ with at least one odd cycle.

\begin{theorem}
\label{thm:five1} Let $G$ be a maximal triangle-free planar graph of
order $11$ with an odd cycle. Then $15\le \varepsilon(G) \le 17$.
\end{theorem}
\smallskip

\noindent{\it Proof.} Consider a planar embedding $G'$ of $G$

Since $G'$ contains odd cycles, from Lemma~\ref{lem:three} we have
$\varepsilon(G')=18-\frac{f_5}{2}\le{17}$. Since $f_4\ge{0}$ we have
$$5 f_5\le 4 f_4 + 5 f_5=2\varepsilon(G)=36 -f_5 .$$
Thus $6 f_5\le 36$ implying that $f_5 \le 6$. Since $f_5>0$ and
even, $f_5=2$, 4 or 6. Thus $\varepsilon(G)=15$, 16 or $17$
according as $f_5=6$, 4 or 2. This completes the proof. \hfill
$\square$

\medskip

We now present the main result of the paper.

\begin{theorem}
\label{thm:six} The smallest order of a triangle-free planar graph
$G$ with $\chi_1(G)=3$ is 11, that is, $f(3,1;\p)=11$.
\end{theorem}

\begin{figure}[H]
 \centering
 \setlength{\unitlength}{1cm}
 \begin{picture}(11,7.5)(0.5,1.5)
  \put(1,5){\circle{0.8}
            \makebox(0,0)[r]{$u$}}
  \put(3,4){%
  \polyline(0,0)(1.5,0)(1.5,3)(0,3)(0,0)}
\put(3.2,6.5){\makebox(0,0)[l]{$u_1$}}
 \put(3.2,5.5){\makebox(0,0)[l]{$u_2$}}
 \put(3.2,4.5){\makebox(0,0)[l]{$u_3$}}
 \put(3.7,6.5){\circle*{0.2}}
 \put(3.7,5.5){\circle*{0.2}}
 \put(3.7,4.5){\circle*{0.2}}
 \put(6.5,4){%
  \polyline(0,0)(1.5,0)(1.5,3)(0,3)(0,0)}
\put(7.7,6.5){\makebox(0,0)[r]{$z_1$}}
 \put(7.7,5.5){\makebox(0,0)[r]{$z_2$}}
 \put(7.7,4.5){\makebox(0,0)[r]{$z_3$}}
 \put(7.2,6.5){\circle*{0.2}}
 \put(7.2,5.5){\circle*{0.2}}
 \put(7.2,4.5){\circle*{0.2}}
 \put(4,2){%
  \polyline(0,0)(3,0)(3,1.5)(0,1.5)(0,0)}
\put(4.5,2.5){\makebox(0,0)[lb]{$u_4$}}
 \put(5.5,2.5){\makebox(0,0)[lb]{$u_5$}}
 \put(6.5,2.5){\makebox(0,0)[lb]{$u_6$}}
 \put(4.5,3.0){\circle*{0.2}}
 \put(5.5,3.0){\circle*{0.2}}
 \put(6.5,3.0){\circle*{0.2}}
 \put(10,5){\circle{0.8}
            \makebox(0,0)[r]{$z$}}
\polyline(3.8,6.5)(7.1,6.5)(3.8,5.5)
 \polyline(7.1,6.5)(3.8,4.5)(7.1,5.5)
\polyline(3.8,4.5)(7.1,4.5)
 \polyline(1.36,5.2)(3.0,6.1)
  \polyline(1.4,5.05)(3.0,5.95)
  \polyline(1.36,4.8)(4,3.2)
  \polyline(1.24,4.7)(4,3.0)
 \polyline(8,6.1)(9.63,5.21)
  \polyline(8,5.95)(9.6,5.06)
  \polyline(7,3.05)(9.6,5.06)
  \polyline(7,2.9)(9.6,4.91)

 \end{picture}
\caption{\label{fig:thmsix} The graph $G^*$ used in
Theorem~\ref{thm:six}.}
\end{figure}
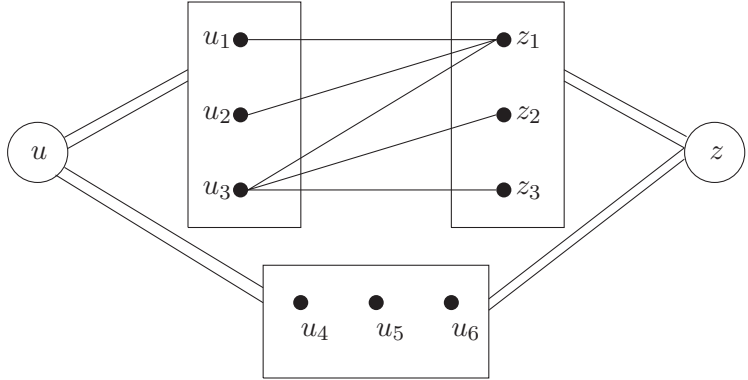

\noindent{\it Proof.} From Theorem~\ref{thm:four} and
Lemma~\ref{lem:two}, it is easy to see that there is no planar
triangle-free graph $G$ on 10 or fewer vertices with $\chi_1(G)=3.$
Hence we have $f(3,1;\p)\ge{11}$. To establish equality, consider
the tfp graph $G_{p1}$ of order 11 shown in Figure~\ref{fig:thmsix}.

\medskip

Suppose $G_{p1}$ is $(2,1)$-colourable and consider a $(2,1)-$
colouring of $G_{p1}$ with colours 1 and 2. Without loss of
generality, we assume that vertex $u$ is assigned colour 1. Apply
Lemma~\ref{lem:one} three times.
\begin{enumerate}
\item[(i)] Applying it first to $G_{p1}[\{u,u_1,u_2,u_3,z_1\}]=K_{2,3}$ we have that
$z_1$ is coloured 1 and hence $u_3$ is assigned colour 2.

\item[(ii)] Now applying  Lemma~\ref{lem:one} to  $G_{p1}[\{u_3,z_1,z_2,z_3,z\}]=K_{2,3}$
it follows that $z$ is assigned colour 2.

\item[(iii)] Finally applying Lemma~\ref{lem:one} to $G_{p1}[\{u,u_4,u_5,u_6,z\}]=K_{2,3}$
it follows that $z$ is assigned colour 1.

\end{enumerate}
But items (ii) and (iii) above are contradictory and so $G_{p1}$ is
not $(2,1)$-colourable. It is easy to see that $G_{p1}$ is
$(3,1)$-colourable (and in fact, $(3,0)$-colourable by
Gr{\"{o}}tzsch theorem). Hence $\chi_1(G_{p1})=3$ and this completes
the proof of Theorem~\ref{thm:six}. \hfill $\square$
\medskip

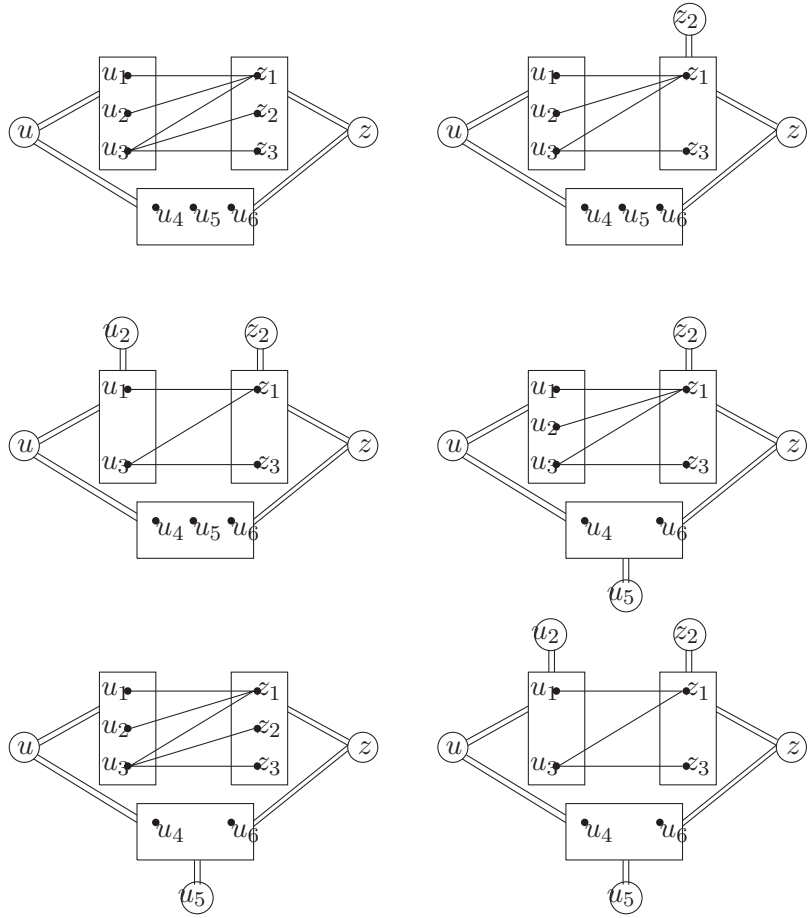
\begin{figure}[H]
 \centering
 \setlength{\unitlength}{0.5cm}
 $$
 \begin{picture}(11,7.5)(0.5,1.5)
  \put(1,5){\circle{0.8}
            \makebox(0,0)[r]{$u$}}
  \put(3,4){%
  \polyline(0,0)(1.5,0)(1.5,3)(0,3)(0,0)}
\put(3.05,6.5){\makebox(0,0)[l]{$u_1$}}
 \put(3.05,5.5){\makebox(0,0)[l]{$u_2$}}
 \put(3.05,4.5){\makebox(0,0)[l]{$u_3$}}
 \put(3.75,6.5){\circle*{0.2}}
 \put(3.75,5.5){\circle*{0.2}}
 \put(3.75,4.5){\circle*{0.2}}
 \put(6.5,4){%
  \polyline(0,0)(1.5,0)(1.5,3)(0,3)(0,0)}
\put(7.8,6.5){\makebox(0,0)[r]{$z_1$}}
 \put(7.8,5.5){\makebox(0,0)[r]{$z_2$}}
 \put(7.8,4.5){\makebox(0,0)[r]{$z_3$}}
 \put(7.2,6.5){\circle*{0.2}}
 \put(7.2,5.5){\circle*{0.2}}
 \put(7.2,4.5){\circle*{0.2}}
 \put(4,2){%
  \polyline(0,0)(3.1,0)(3.1,1.5)(0,1.5)(0,0)}
\put(4.5,2.5){\makebox(0,0)[lb]{$u_4$}}
 \put(5.5,2.5){\makebox(0,0)[lb]{$u_5$}}  
 \put(6.5,2.5){\makebox(0,0)[lb]{$u_6$}}
 \put(4.5,3.0){\circle*{0.2}}
 \put(5.5,3.0){\circle*{0.2}}
 \put(6.5,3.0){\circle*{0.2}}
 \put(10,5){\circle{0.8}
            \makebox(0,0)[r]{$z$}}
\polyline(3.8,6.5)(7.1,6.5)(3.8,5.5)
 \polyline(7.1,6.5)(3.8,4.5)(7.1,5.5)
\polyline(3.8,4.5)(7.1,4.5)
 \polyline(1.36,5.2)(3.0,6.1)
  \polyline(1.4,5.05)(3.0,5.95)
  \polyline(1.36,4.8)(4,3.2)
  \polyline(1.24,4.7)(4,3.0)
 \polyline(8,6.1)(9.63,5.21)
  \polyline(8,5.95)(9.6,5.06)
  \polyline(7.1,3.05)(9.6,5.06)
  \polyline(7.1,2.9)(9.6,4.91)

 \end{picture}
 \hspace{0.2cm}
 \begin{picture}(11,7.5)(0.5,1.5)
  \put(1,5){\circle{0.8}
            \makebox(0,0)[r]{$u$}}
  \put(3,4){%
  \polyline(0,0)(1.5,0)(1.5,3)(0,3)(0,0)}
\put(3.05,6.5){\makebox(0,0)[l]{$u_1$}}
 \put(3.05,5.5){\makebox(0,0)[l]{$u_2$}}
 \put(3.05,4.5){\makebox(0,0)[l]{$u_3$}}
 \put(3.75,6.5){\circle*{0.2}}
 \put(3.75,5.5){\circle*{0.2}}
 \put(3.75,4.5){\circle*{0.2}}
 \put(6.5,4){%
  \polyline(0,0)(1.5,0)(1.5,3)(0,3)(0,0)}
\put(7.85,6.5){\makebox(0,0)[r]{$z_1$}}
 \put(7.85,4.5){\makebox(0,0)[r]{$z_3$}}
 \put(7.2,6.5){\circle*{0.2}}
 \put(7.2,4.5){\circle*{0.2}}
  \put(7.3,8){\circle{0.85}
            \makebox(0,0)[r]{$z_2$}}
 \put(4,2){%
  \polyline(0,0)(3.1,0)(3.1,1.5)(0,1.5)(0,0)}
\put(4.5,2.5){\makebox(0,0)[lb]{$u_4$}}
 \put(5.5,2.5){\makebox(0,0)[lb]{$u_5$}}  
 \put(6.5,2.5){\makebox(0,0)[lb]{$u_6$}}
 \put(4.5,3.0){\circle*{0.2}}
 \put(5.5,3.0){\circle*{0.2}}
 \put(6.5,3.0){\circle*{0.2}}
 \put(10,5){\circle{0.8}
            \makebox(0,0)[r]{$z$}}
\polyline(3.8,6.5)(7.1,6.5)(3.8,5.5)
\polyline(3.8,4.5)(7.1,6.5)
\polyline(3.8,4.5)(7.1,4.5)
 \polyline(1.36,5.2)(3.0,6.1)
  \polyline(1.4,5.05)(3.0,5.95)
  \polyline(1.36,4.8)(4,3.2)
  \polyline(1.24,4.7)(4,3.0)
 \polyline(8,6.1)(9.63,5.21)
  \polyline(8,5.95)(9.6,5.06)
  \polyline(7.1,3.05)(9.6,5.06)
  \polyline(7.1,2.9)(9.6,4.91)
\polyline(7.35,7)(7.35,7.6) \polyline(7.2,7)(7.2,7.6)

 \end{picture}
 $$
 $$
  \begin{picture}(11,7.5)(0.5,1.5)
  \put(1,5){\circle{0.8}
            \makebox(0,0)[r]{$u$}}
  \put(3,4){%
  \polyline(0,0)(1.5,0)(1.5,3)(0,3)(0,0)}
\put(3.05,6.5){\makebox(0,0)[l]{$u_1$}}
 \put(3.05,4.5){\makebox(0,0)[l]{$u_3$}}
 \put(3.75,6.5){\circle*{0.2}}
 \put(3.75,4.5){\circle*{0.2}}
   \put(3.6,8){\circle{0.85}
            \makebox(0,0)[r]{$u_2$}}
 \put(6.5,4){%
  \polyline(0,0)(1.5,0)(1.5,3)(0,3)(0,0)}
\put(7.85,6.5){\makebox(0,0)[r]{$z_1$}}
 \put(7.85,4.5){\makebox(0,0)[r]{$z_3$}}
 \put(7.2,6.5){\circle*{0.2}}
 \put(7.2,4.5){\circle*{0.2}}
  \put(7.3,8){\circle{0.85}
            \makebox(0,0)[r]{$z_2$}}
 \put(4,2){%
  \polyline(0,0)(3.1,0)(3.1,1.5)(0,1.5)(0,0)}
\put(4.5,2.5){\makebox(0,0)[lb]{$u_4$}}
 \put(5.5,2.5){\makebox(0,0)[lb]{$u_5$}}  
 \put(6.5,2.5){\makebox(0,0)[lb]{$u_6$}}
 \put(4.5,3.0){\circle*{0.2}}
 \put(5.5,3.0){\circle*{0.2}}
 \put(6.5,3.0){\circle*{0.2}}
 \put(10,5){\circle{0.8}
            \makebox(0,0)[r]{$z$}}
\polyline(3.8,6.5)(7.1,6.5)
\polyline(3.8,4.5)(7.1,6.5)
\polyline(3.8,4.5)(7.1,4.5)
 \polyline(1.36,5.2)(3.0,6.1)
  \polyline(1.4,5.05)(3.0,5.95)
  \polyline(1.36,4.8)(4,3.2)
  \polyline(1.24,4.7)(4,3.0)
 \polyline(8,6.1)(9.63,5.21)
  \polyline(8,5.95)(9.6,5.06)
  \polyline(7.1,3.05)(9.6,5.06)
  \polyline(7.1,2.9)(9.6,4.91)
\polyline(7.35,7)(7.35,7.6) \polyline(7.2,7)(7.2,7.6)
\polyline(3.55,7)(3.55,7.6) \polyline(3.7,7)(3.7,7.6)

 \end{picture}
 \hspace{0.2cm}
  \begin{picture}(11,7.5)(0.5,1.5)
  \put(1,5){\circle{0.8}
            \makebox(0,0)[r]{$u$}}
  \put(3,4){%
  \polyline(0,0)(1.5,0)(1.5,3)(0,3)(0,0)}
\put(3.05,6.5){\makebox(0,0)[l]{$u_1$}}
 \put(3.05,5.5){\makebox(0,0)[l]{$u_2$}}
 \put(3.05,4.5){\makebox(0,0)[l]{$u_3$}}
 \put(3.75,6.5){\circle*{0.2}}
 \put(3.75,5.5){\circle*{0.2}}
 \put(3.75,4.5){\circle*{0.2}}
 \put(6.5,4){%
  \polyline(0,0)(1.5,0)(1.5,3)(0,3)(0,0)}
\put(7.85,6.5){\makebox(0,0)[r]{$z_1$}}
 \put(7.85,4.5){\makebox(0,0)[r]{$z_3$}}
 \put(7.2,6.5){\circle*{0.2}}
 \put(7.2,4.5){\circle*{0.2}}
  \put(7.3,8){\circle{0.85}
            \makebox(0,0)[r]{$z_2$}}
 \put(4,2){%
  \polyline(0,0)(3.1,0)(3.1,1.5)(0,1.5)(0,0)}
\put(4.5,2.5){\makebox(0,0)[lb]{$u_4$}}
 \put(6.5,2.5){\makebox(0,0)[lb]{$u_6$}}
 \put(4.5,3.0){\circle*{0.2}}
 \put(6.5,3.0){\circle*{0.2}}
  \put(5.6,1){\circle{0.85}
            \makebox(0,0)[r]{$u_5$}}
\polyline(5.525,1.35)(5.525,2) \polyline(5.675,1.35)(5.675,2)
 \put(10,5){\circle{0.8}
            \makebox(0,0)[r]{$z$}}
\polyline(3.8,6.5)(7.1,6.5)(3.8,5.5)
\polyline(3.8,4.5)(7.1,6.5)
\polyline(3.8,4.5)(7.1,4.5)
 \polyline(1.36,5.2)(3.0,6.1)
  \polyline(1.4,5.05)(3.0,5.95)
  \polyline(1.36,4.8)(4,3.2)
  \polyline(1.24,4.7)(4,3.0)
 \polyline(8,6.1)(9.63,5.21)
  \polyline(8,5.95)(9.6,5.06)
  \polyline(7.1,3.05)(9.6,5.06)
  \polyline(7.1,2.9)(9.6,4.91)
\polyline(7.35,7)(7.35,7.6) \polyline(7.2,7)(7.2,7.6)

 \end{picture}
 $$
  $$
 \begin{picture}(11,7.5)(0.5,1.5)
  \put(1,5){\circle{0.8}
            \makebox(0,0)[r]{$u$}}
  \put(3,4){%
  \polyline(0,0)(1.5,0)(1.5,3)(0,3)(0,0)}
\put(3.05,6.5){\makebox(0,0)[l]{$u_1$}}
 \put(3.05,5.5){\makebox(0,0)[l]{$u_2$}}
 \put(3.05,4.5){\makebox(0,0)[l]{$u_3$}}
 \put(3.75,6.5){\circle*{0.2}}
 \put(3.75,5.5){\circle*{0.2}}
 \put(3.75,4.5){\circle*{0.2}}
 \put(6.5,4){%
  \polyline(0,0)(1.5,0)(1.5,3)(0,3)(0,0)}
\put(7.85,6.5){\makebox(0,0)[r]{$z_1$}}
 \put(7.85,5.5){\makebox(0,0)[r]{$z_2$}}
 \put(7.85,4.5){\makebox(0,0)[r]{$z_3$}}
 \put(7.2,6.5){\circle*{0.2}}
 \put(7.2,5.5){\circle*{0.2}}
 \put(7.2,4.5){\circle*{0.2}}
 \put(4,2){%
  \polyline(0,0)(3.1,0)(3.1,1.5)(0,1.5)(0,0)}
\put(4.5,2.5){\makebox(0,0)[lb]{$u_4$}}
 \put(6.5,2.5){\makebox(0,0)[lb]{$u_6$}}
 \put(4.5,3.0){\circle*{0.2}}
 \put(6.5,3.0){\circle*{0.2}}
   \put(5.6,1){\circle{0.85}
            \makebox(0,0)[r]{$u_5$}}
\polyline(5.525,1.35)(5.525,2) \polyline(5.675,1.35)(5.675,2)
 \put(10,5){\circle{0.8}
            \makebox(0,0)[r]{$z$}}
\polyline(3.8,6.5)(7.1,6.5)(3.8,5.5)
 \polyline(7.1,6.5)(3.8,4.5)(7.1,5.5)
\polyline(3.8,4.5)(7.1,4.5)
 \polyline(1.36,5.2)(3.0,6.1)
  \polyline(1.4,5.05)(3.0,5.95)
  \polyline(1.36,4.8)(4,3.2)
  \polyline(1.24,4.7)(4,3.0)
 \polyline(8,6.1)(9.63,5.21)
  \polyline(8,5.95)(9.6,5.06)
  \polyline(7.1,3.05)(9.6,5.06)
  \polyline(7.1,2.9)(9.6,4.91)

 \end{picture}
 \hspace{0.2cm}
   \begin{picture}(11,7.5)(0.5,1.5)
  \put(1,5){\circle{0.8}
            \makebox(0,0)[r]{$u$}}
  \put(3,4){%
  \polyline(0,0)(1.5,0)(1.5,3)(0,3)(0,0)}
\put(3.05,6.5){\makebox(0,0)[l]{$u_1$}}
 \put(3.05,4.5){\makebox(0,0)[l]{$u_3$}}
 \put(3.75,6.5){\circle*{0.2}}
 \put(3.75,4.5){\circle*{0.2}}
   \put(3.6,8){\circle{0.85}
            \makebox(0,0)[r]{$u_2$}}
 \put(6.5,4){%
  \polyline(0,0)(1.5,0)(1.5,3)(0,3)(0,0)}
\put(7.85,6.5){\makebox(0,0)[r]{$z_1$}}
 \put(7.85,4.5){\makebox(0,0)[r]{$z_3$}}
 \put(7.2,6.5){\circle*{0.2}}
 \put(7.2,4.5){\circle*{0.2}}
  \put(7.3,8){\circle{0.85}
            \makebox(0,0)[r]{$z_2$}}
 \put(4,2){%
  \polyline(0,0)(3.1,0)(3.1,1.5)(0,1.5)(0,0)}
\put(4.5,2.5){\makebox(0,0)[lb]{$u_4$}}
 \put(6.5,2.5){\makebox(0,0)[lb]{$u_6$}}
 \put(4.5,3.0){\circle*{0.2}}
 \put(6.5,3.0){\circle*{0.2}}
   \put(5.6,1){\circle{0.85}
            \makebox(0,0)[r]{$u_5$}}
\polyline(5.525,1.35)(5.525,2) \polyline(5.675,1.35)(5.675,2)
 \put(10,5){\circle{0.8}
            \makebox(0,0)[r]{$z$}}
\polyline(3.8,6.5)(7.1,6.5)
\polyline(3.8,4.5)(7.1,6.5)
\polyline(3.8,4.5)(7.1,4.5)
 \polyline(1.36,5.2)(3.0,6.1)
  \polyline(1.4,5.05)(3.0,5.95)
  \polyline(1.36,4.8)(4,3.2)
  \polyline(1.24,4.7)(4,3.0)
 \polyline(8,6.1)(9.63,5.21)
  \polyline(8,5.95)(9.6,5.06)
  \polyline(7.1,3.05)(9.6,5.06)
  \polyline(7.1,2.9)(9.6,4.91)
\polyline(7.35,7)(7.35,7.6) \polyline(7.2,7)(7.2,7.6)
\polyline(3.55,7)(3.55,7.6) \polyline(3.7,7)(3.7,7.6)

 \end{picture}
 $$
\caption{\label{fig:thmseven} The 6 graphs $G_{pi}$ for $1\le
i\le6$.}
\end{figure}

\medskip

The graph $G_{p6}-(u_{3},z_{1})$ is isomorphic to the ``cross-over"
graph of Figure 3 of ~\cite{FJLS}.

\begin{lemma}
\label{lem:five}
 Let $G$ be a tfp graph of order 11 with $\delta (G)= 2$. Let $e$ be
 an edge incident with a vertex of degree 2. Then $\chi_{1}(G-e)= 2.$
\end{lemma}

\noindent{\it Proof.} Suppose that $e=(x,y)$ and $d_{G}(x)=2.$ Let
$x_1$ be the other neighbour of $x.$ The graph $G-x$ is tfp graph of
order 10 and hence $\chi_{1}(G-x)=2.$ Consider a $(2,1)-$colouring
of $G-x$ using colours 1 and 2. Without loss of generality assume
that the vertex $x_1$ is assigned colour 1. Now assign colour 2 to
vertex $x.$ This provides a $(2,1)-$colouring of $G-e$ and hence
$\chi_{1}(G-e)=2$. This completes the proof of Lemma~\ref{lem:five}.
\hfill $\square$

\medskip

The result of the next lemma was evident from the computations, but
here we give a simple direct proof.

\begin{lemma}
\label{lem:six}
 The tfp graph $G_{p6}$ of order 11 given in Figure~\ref{fig:thmseven}
 is $(3,1)-$edge-critical.
\end{lemma}

\noindent{\it Proof.} Let $e=(x,y)$ be an edge of $G$. By
Lemma~\ref{lem:five}, we can assume that $d_{G}(x)> 2$ and
$d_{G}(y)> 2$. There are precisely 11 such edges in $G_{p6}$ and
they are one of five types. For a representative edge $e=(x,y)$ of
each type the table below describes a $(2,1)-$colouring of $G-e$.
The $(2,1)-$colouring is defined by the partition of $V(G)$ into
$1-$independent sets $X(e)$ and $V(G)-X(e)$.

\begin{center}
\begin{tabular}{|l|l|l|}
\hline
$e$& $X(e)$& other edges of \\
 & & same type\\
 \hline
$(u,u_{4})$&$\{u,z_1,z_3,u_2,u_4,u_6 \}$&
$(u,u_6),(z,u_4),(z,u_6)$\\
$(u,u_1)$&$\{u,z,u_1,u_3,u_5, z_2 \}$& $(z,z_3)$\\
$(u_1,z_1)$&$\{u,z,u_3,u_5, z_2 \}$& $(u_3,z_3)$\\
$(u,u_3)$&$\{u,z,u_1,u_3,u_5, z_2 \}$& $(z,z_1)$\\
$(u_3,z_1)$&$\{u,z,u_2,u_5, z_2 \}$& \\
\hline
 \end{tabular}
\end{center}
This establishes the Lemma~\ref{lem:six}. \hfill $\square$

\begin{theorem}
\label{thm:eight} The tfp graphs $G_{pi}$ for $1\le i \le 6$ are
$(3,1)-$edge-critical.
\end{theorem}

\noindent{\it Proof.} The Lemma~\ref{lem:six} establishes the result
for graph $G_{p6}$. The results for the other $G_{pi}$ are proved in
a similar manner. These cases are simpler as for each of them the
number of edges free of degree $2$ is smaller than it was for the
graph $G_{p6}$. More precisely there are just 3, 5, 5 , 9 and 7 such
edges respectively for the graphs $G_{pi}$, $1\le i \le{5}$.\hfill
$\square$

\end{document}